\documentclass[12pt]{amsart}
\usepackage{amsmath,amscd,amssymb,amsfonts}
\newcommand{\ver}{Oct. 25, 2006, v.1}
\newcommand{\ssbull}{\raise.2ex\hbox{${\scriptstyle\bullet}$}}
\newcommand{\simto}{\buildrel\sim\over\longrightarrow}
\newcommand{\msum}{\hbox{$\sum$}}
\newcommand{\mprod}{\hbox{$\prod$}}
\newcommand{\moplus}{\hbox{$\bigoplus$}}
\newcommand{\mcup}{\hbox{$\bigcup$}}

\newcommand{\bs}{{\mathbf s}}
\newcommand{\bC}{{\mathbf C}}
\newcommand{\bD}{{\mathbf D}}
\newcommand{\bN}{{\mathbf N}}
\newcommand{\bP}{{\mathbf P}}
\newcommand{\bQ}{{\mathbf Q}}
\newcommand{\bR}{{\mathbf R}}
\newcommand{\bZ}{{\mathbf Z}}
\newcommand{\cA}{{\mathcal A}}
\newcommand{\cD}{{\mathcal D}}
\newcommand{\cE}{{\mathcal E}}
\newcommand{\cH}{{\mathcal H}}
\newcommand{\cI}{{\mathcal I}}
\newcommand{\cJ}{{\mathcal J}}
\newcommand{\cL}{{\mathcal L}}
\newcommand{\cO}{{\mathcal O}}
\newcommand{\fa}{{\mathfrak a}}
\newcommand{\tR}{\widetilde{R}}
\newcommand{\tb}{\widetilde{b}}

\newcommand{\tF}{\widetilde{F}}
\newcommand{\tm}{\widetilde{m}}
\newcommand{\tD}{\widetilde{D}}
\newcommand{\tM}{\widetilde{M}}
\newcommand{\tP}{\widetilde{P}}
\newcommand{\tH}{\widetilde{H}}
\newcommand{\tU}{\widetilde{U}}

\newcommand{\tX}{\widetilde{X}}

\newcommand{\tal}{\widetilde{\alpha}}
\newcommand{\Gr}{\hbox{\rm Gr}}
\newcommand{\DR}{\hbox{\rm DR}}
\newcommand{\cHom}{\hbox{${\mathcal H}om$}}
\newcommand{\Ker}{\hbox{\rm Ker}}
\newcommand{\Perv}{\hbox{\rm Perv}}
\newcommand{\cSpec}{\hbox{${\mathcal S}pec$}}
\newcommand{\Sp}{\hbox{\rm Sp}}
\newcommand{\Sym}{\hbox{\rm Sym}}
\newcommand{\Sing}{\hbox{\rm Sing}}
\newcommand{\Supp}{\hbox{\rm Supp}}
\newcommand{\CV}{\hbox{\rm CV}}
\newcommand{\wdg}{\hbox{$\wedge$}}
\newcommand{\rd}{\partial}
\newcommand{\finv}{\hbox{$\frac{1}{f}$}}
\newcommand{\Mustata}{Musta\c{t}\v{a}}

\begin{document}
\title{Introduction to a theory of $b$-functions}
\author{Morihiko Saito}
\address{RIMS Kyoto University, Kyoto 606-8502 Japan}
\email{msaito@kurims.kyoto-u.ac.jp}
\date{\ver}
%\begin{abstract}
%\end{abstract}
\maketitle

\noindent
We give an introduction to a theory of $b$-functions, i.e.
Bernstein-Sato polynomials. After reviewing some facts from
$D$-modules, we introduce $b$-functions including the one for
arbitrary ideals of the structure sheaf. We explain the relation with
singularities, multiplier ideals, etc., and calculate the $b$-functions
of monomial ideals and also of hyperplane arrangements in certain cases.

\bigskip\bigskip
\centerline{\bf 1. D-modules.}

\bigskip\noindent
{\bf 1.1.}
Let $X$ be a complex manifold or a smooth algebraic variety over $\bC$.
Let $\cD_X$ be the ring of partial differential operators.
A local section of $\cD_X$ is written as
$$
\msum_{\nu\in\bN^n} \,a_{\nu}\rd_1^{\nu_1}\cdots
\rd_n^{\nu_n}\in\cD_X\quad\text{with}\,\,a_{\nu}\in\cO_X,
$$
where $\rd_i=\rd/\rd x_i$ with $(x_1,\dots,x_n)$ a local coordinate
system.

Let $F$ be the filtration by the order of operators i.e.
$$
F_p\cD_X=\big\{\msum_{|\nu|\le p} \,a_{\nu}\rd_1^{\nu_1}
\cdots\rd_n^{\nu_n}\big\},
$$
where $|\nu|=\sum_i\nu_i$.
Let $\xi_i=\Gr^F_1\rd_i\in\Gr^F_1\cD_X$.
Then
$$
\aligned
&\Gr^F\cD_X:=\moplus_p\Gr^F_p\cD_X=\moplus_p\Sym^p\Theta_X
\,(=\cO_X[\xi_1,\dots,\xi_n]\,\,\,\text{locally)},
\\
&\cSpec_X\Gr^F\cD_X=T^*X.
\endaligned
\leqno(1.1.1)
$$

\medskip\noindent
{\bf 1.2~Definition.}
We say that a left $\cD_X$-module $M$ is {\it coherent}
if it has locally a finite presentation
$$
\moplus \cD_X\to \moplus\cD_X\to M\to 0.
$$

\medskip\noindent
{\bf 1.3.~Remark.}
A left $\cD_X$-module $M$ is coherent if and only if it is
quasi-coherent over $\cO_X$ and locally finitely generated over
$\cD_X$.
(It is known that $\Gr^F\cD_X$ is a noetherian ring, i.e.
an increasing sequence of locally finitely generated
$\Gr^F\cD_X$-submodules of a coherent $\Gr^F\cD_X$-module
is locally stationary.)

\medskip\noindent
{\bf 1.4.~Definition.}
A filtration $F$ on a left $\cD_X$-module $M$ is {\it good}
if $(M,F)$ is a coherent filtered $\cD_X$-module, i.e.
if $F_p\cD_XF_qM\subset M_{p+q}$ and
$\Gr^FM:=\moplus_p\Gr_p^FM$ is coherent over $\Gr^F\cD_X$.

\medskip\noindent
{\bf 1.5.~Remark.} A left $\cD_X$-module $M$ is coherent if and only if
it has a good filtration locally.

\medskip\noindent
{\bf 1.6.~Characteristic varieties.}
For a coherent left $\cD_X$-module $M$, we define the characteristic
variety $\CV(M)$ by
$$
\CV(M) = \Supp\,\Gr^FM\subset T^*M,
\leqno(1.6.1)
$$
taking locally a good filtration $F$ of $M$.

\medskip\noindent
{\bf 1.7.~Remark.}
The above definition is independent of the choice of $F$.
If $M=\cD_X/\cI$ for a coherent left ideal $\cI$ of $\cD_X$,
take $P_i\in F_{k_i}\cI$ such that the $\rho_i:=\Gr^F_{k_i}P_i$
generate $\Gr^F\cI$ over $\Gr^F\cD_X$.
Then $\CV(M)$ is defined by the $\rho_i\in\cO_X[\xi_1,\dots,\xi_n]$.

\medskip\noindent
{\bf 1.8.~Theorem} (Sato, Kawai, Kashiwara [39], Bernstein [2]).
{\it We have the inequality $\dim \CV(M)\ge\dim X$.
{\rm (}More precisely, $\CV(M)$ is involutive, see {\rm [39].)}
}

\medskip\noindent
{\bf 1.9.~Definition.}
We say that a left $\cD_X$-module
$M$ is {\it holonomic} if it is coherent and
$\dim \CV(M)=\dim X$.

\bigskip\bigskip
\centerline{\bf 2. De Rham functor.}

\bigskip\noindent
{\bf 2.1.~Definition.}
For a left $\cD_X$-module $M$,
we define the de Rham functor $\DR(M)$ by
$$
M\to\Omega_X^1\otimes_{\cO_X}M\to\cdots\to\Omega_X^{\dim X}
\otimes_{\cO_X}M,
\leqno(2.1.1)
$$
where the last term is put at the degree 0.
In the algebraic case, we use analytic sheaves or replace
$M$ with the associated analytic sheaf $M^{\rm an}:=M\otimes_{\cO_X}
\cO_{X^{\rm an}}$ in case $M$ is algebraic (i.e. $M$ is an
$\cO_X$-module with $\cO_X$ algebraic).

\medskip\noindent
{\bf 2.2.~Perverse sheaves.}
Let
$D_c^b(X,\bC)$ be the derived category of bounded complexes
of $\bC_X$-modules $K$ with $\cH^jK$ constructible.
(In the algebraic case we use analytic topology for the sheaves
although we use Zariski topology for constructibility.)
Then the category of perverse sheaves
$\Perv(X,\bC)$ is a full subcategory of $D_c^b(X,\bC)$ consisting
of $K$ such that
$$
\dim \Supp\,\cH^{-j}K\le j,
\quad\dim \Supp\,\cH^{-j}\bD K\le j,
\leqno(2.2.1)
$$
where $\bD K:=\bR\cHom(K,\bC[2\dim X])$ is the dual of $K$, and
$\cH^jK$ is the $j$-th cohomology sheaf of $K$.

\medskip\noindent
{\bf 2.3.~Theorem} (Beilinson, Bernstein, Deligne [1]). {\it
$\Perv(X,\bC)$ is an abelian category.
}

\medskip\noindent
{\bf 2.4.~Theorem} (Kashiwara). {\it
If $M$ is holonomic, then
$\DR(M)$ is a perverse sheaf.
}

\medskip\noindent
{\it Outline of proof.}
By Kashiwara [19], we have $\DR(M)\in D^b_c(X,\bC)$, and
the first condition of (2.2.1) is verified.
Then the assertion follows from the commutativity of the dual
$\bD$ and the de Rham functor $\DR$.

\medskip\noindent
{\bf 2.5.~Example.} $\DR(\cO_X)=\bC_X[\dim X]$.

\medskip\noindent
{\bf 2.6.~Direct images.}
For a closed immersion $i:X\to Y$ such that
$X$ is defined by $x_i=0$ in $Y$ for $1\le i\le r$, define the direct
image of left $\cD_X$-modules $M$ by
$$
i_+M:=M[\rd_1,\dots,\rd_r].
$$
(Globally there is a twist by a line bundle.)
For a projection $p:X\times Y\to Y$, define
$$
p_+M=\bR p_*\DR_X(M).
$$
In general, $f_+=p_+i_+$ using $f=pi$
with $i$ graph embedding.
See [4] for details.

\medskip\noindent
{\bf 2.7.~Regular holonomic D-modules.}
Let $M$ be a holonomic $\cD_X$-module with support $Z$,
and $U$ be a Zariski-open of $Z$ such that
$\DR(M)|_U$ is a local system up to a shift.
Then $M$ is {\it regular} if and only if there exists locally a
divisor $D$ on $X$ containing $Z\setminus U$ and such that $M(*D)$
is the direct image of
a regular holonomic $\cD$-module `of Deligne-type' (see [11])
on a desingularization of $(Z,Z\cap D)$,
and Ker($M\to M(*D)$) is regular holonomic
(by induction on $\dim \Supp\, M$).

Note that the category $M_{rh}(\cD_X)$ of regular holonomic
$\cD_X$-modules is stable by subquotients and extensions
in the category $M_{h}(\cD_X)$ of holonomic $\cD_X$-modules.

\medskip\noindent
{\bf 2.8.~Theorem} (Kashiwara-Kawai [24], [22], Mebkhout [28]).

\noindent
(i) The structure sheaf $\cO_X$ is regular holonomic.

\noindent
(ii) The functor DR induces an equivalence of categories
$$
\DR:M_{rh}(\cD_X)\simto\Perv(X,\bC).
\leqno(2.8.1)
$$

\noindent
(See [4] for the algebraic case.)

\bigskip\bigskip
\centerline{\bf 3. $b$-Functions.}

\bigskip\noindent
{\bf 3.1.~ Definition.}
Let
$f$ be a holomorphic function on $X$, or $f\in \Gamma(X,\cO_X)$
in the algebraic case.
Then we have
$$
\cD_X[s]f^s\subset \cO_X[{\finv}][s]f^s\quad\text{where}\,\,
\rd_i f^s=s(\rd_if)f^{s-1},
$$
and $b_f(s)$ is the monic polynomial of the least degree satisfying
$$
b_f(s)f^s=P(x,\rd,s)f^{s+1}\quad\text{in}\,\,
\cO_X[{\finv}][s]f^s,
$$
with $P(x,\rd,s)\in\cD_X[s]$. Locally, it is  the minimal
polynomial of the action of $s$ on
$$
\cD_X[s]f^s/\cD_X[s]f^{s+1}.
$$
We define $b_{f,x}(s)$ replacing $\cD_X$ with $\cD_{X,x}$.

\medskip\noindent
{\bf 3.2.~Theorem} (Sato [38], Bernstein [2], Bjork [3]). {\it
The $b$-function exists
at least locally, and exists globally in the case $X$ affine
variety with $f$ algebraic.
}

\medskip\noindent
{\bf 3.3.~Observation.}
Let $i_f:X\to \tX:=X\times\bC$ be the graph embedding.
Then there are canonical isomorphisms
$$
\tM:=i_{f+}\cO_X=\cO_X[\rd_t]\delta(f-t)
=\cO_{X\times\bC}\bigl[\hbox{$\frac{1}{f-t}$}\bigr]/\cO_{X\times\bC},
\leqno(3.3.1)
$$
where the action of $\rd_i$ on $\delta(f-t)\,(=\frac{1}{f-t})$ is given by
$$
\rd_i\delta(f-t)=-(\rd_if)\rd_t\delta(f-t).
\leqno(3.3.2)
$$
Moreover,
$f^s$ is canonically identified with $\delta(f-t)$ setting
$s = -\rd_tt$, and we have a canonical isomorphism as $\cD_X[s]$-modules
$$
\cD_X[s]f^s=\cD_X[s]\delta(f-t).
\leqno(3.3.3)
$$

\medskip\noindent
{\bf 3.4.~V-filtration.}
We say that
$V$ is a filtration of Kashiwara-Malgrange if $V$ is
exhaustive, separated, and satisfies for any
$\alpha\in\bQ$:

\smallskip
(i) $V^{\alpha}\tM$ is a coherent $\cD_X[s]$-submodule of $\tM$.

(ii) $tV^{\alpha}\tM\subset V^{\alpha+1}\tM$ and $=$ holds for
$\alpha\gg 0$.

(iii) $\rd_tV^{\alpha}\tM\subset V^{\alpha-1}\tM$.

(iv) $\rd_tt-\alpha$ is nilpotent on $\Gr_V^{\alpha}\tM$.

\smallskip\noindent
If it exists, it is unique.

\medskip\noindent
{\bf 3.5.~Relation with the $b$-function.} If $X$ is affine or
Stein and relatively compact, then the multiplicity of a root
$\alpha$ of $b_f(s)$ is given by the minimal polynomial of
$s-\alpha$ on
$$
\Gr_V^{\alpha}(\cD_X[s]f^s/\cD_X[s]f^{s+1}),
\leqno(3.5.1)
$$
using $\cD_X[s]f^s=\cD_X[s]\delta(f-t)$ with $s=-\rd_tt$.

Note that
$V^{\alpha}\tM$ and $\cD_X[s]f^{s+i}$ are `lattices' of $\tM$, i.e.
$$
V^{\alpha}\tM\subset\cD_X[s]f^{s+i}\subset V^{\beta}\tM\quad
\text{for}\,\,\alpha\gg i\gg \beta,
\leqno(3.5.2)
$$
and $V^{\alpha}\tM$ is an analogue of the Deligne extension
with eigenvalues in $[\alpha,\alpha+1)$.
The  existence of $V$ is equivalent to the existence of
$b_f(s)$ locally.

\medskip\noindent
{\bf 3.6.~Theorem} (Kashiwara [21], [23], Malgrange [27]). {\it
The filtration $V$ exists on $\tM:=i_{f+}M$ for any holonomic
$\cD_X$-module $M$.
}

\medskip\noindent
{\bf 3.7.~Remarks.}
(i) There are many ways to prove this theorem,
since it is essentially equivalent to the existence of
the $b$-function (in a generalized sense). One way is to
use a resolution of singularities and reduce to the case
where $\CV(M)$ has normal crossings, if $M$ is regular.

(ii) The filtration $V$ is indexed by $\bQ$ if $M$ is quasi-unipotent.

\medskip\noindent
{\bf 3.8.~Relation with vanishing cycle functors.}
Let $\rho:X_t\to X_0$ be a `good' retraction
(using a resolution of singularities of $(X,X_0))$,
where $X_t=f^{-1}(t)$ with $t\ne 0$ sufficiently near 0.
Then we have canonical isomorphisms
$$
\psi_f\bC_X=\bR\rho_*\bC_{X_t},\quad
\varphi_f\bC_X=\psi_f\bC_X/\bC_{X_0},
\leqno(3.8.1)
$$
where $\psi_f\bC_X,\varphi_f\bC_X$ are nearby and vanishing cycle
sheaves, see [13].

Let $F_x$ denote the Milnor fiber around $x\in X_0$. Then
$$
(\cH^j\psi_f\bC_X)_x=H^j(F_x,\bC),\quad
(\cH^j\varphi_f\bC_X)_x=\tH^j(F_x,\bC).
\leqno(3.8.2)
$$

For a $\cD_X$-module $M$ admitting the V-filtration on
$\tM=i_{*+}M$, we define $\cD_X$-modules
$$
\psi_f M=\moplus_{0<\alpha\le 1}\Gr_V^{\alpha}\tM,\quad
\varphi_f M=\moplus_{0\le\alpha<1}\Gr_V^{\alpha}\tM.
\leqno(3.8.3)
$$

\medskip\noindent
{\bf 3.9.~Theorem} (Kashiwara [23], Malgrange [27]). {\it
For a regular holonomic $\cD_X$-module $M$, we have
canonical isomorphisms
$$
\aligned
\DR_X\psi_f(M)&=\psi_f\DR_{X}(M)[-1],
\\
\DR_X\varphi_f(M)&=\varphi_f\DR_{X}(M)[-1],
\endaligned
\leqno(3.9.1)
$$
and $\exp(-2\pi i\rd_tt)$ on the left-hand side corresponds to
the monodromy $T$ on the right-hand side.
}

\medskip\noindent
{\bf 3.10.~Definition.} Let

\smallskip
$R_f= \{$roots of $b_f(-s)\}$,

$\alpha_f=\min R_f$,

$m_{\alpha}$ : the multiplicity of $\alpha \in R_f$.

\smallskip\noindent
(Similarly for $R_{f,x}$, etc. for $b_{f,x}(s)$.)

\medskip\noindent
{\bf 3.11.~Theorem} (Kashiwara [20]). $R_f\subset \bQ_{>0}$.

\medskip
(This is proved by using a resolution of singularities.)

\medskip\noindent
{\bf 3.12.~Theorem} (Kashiwara [23], Malgrange [27]). 

\smallskip
(i) {\it $e^{-2\pi iR_f}=\{$the eigenvalues of $T$ on $H^j(F_x,\bC)$
for $x\in X_0, j\in \bZ\}$,}

\smallskip
(ii) {\it $m_{\alpha}\le\min\{i\mid N^i\psi_{f,\lambda}\bC_X= 0\}$
with $\lambda=e^{-2\pi i\alpha}$,

\smallskip
where $\psi_{f,\lambda}=\Ker(T_s-\lambda)\subset \psi_f$,
$N=\log T_u$ with
$T=T_sT_u$.}

\medskip\noindent
(This is a corollary of the above Theorem (3.9) of Kashiwara
and Malgrange.)

\bigskip\bigskip
\centerline{\bf 4. Relation with other invariants.}

\bigskip\noindent
{\bf 4.1.~Microlocal $b$-function.}
We define $\tR_f,\tm_{\alpha},\tal_f$ with $b_f(s)$ replaced
by the
{\it microlocal} (or reduced) $b$-function
$$
\tb_f(s):=b_f(s)/(s+1).
\leqno(4.1.1)
$$
This $\tb_f(s)$ coincides with the monic polynomial of the least
degree satisfying
$$
\tb_f(s)\delta(f-t)=\tP\rd_t^{-1}\delta(f-t)\quad\text{with}\,\,\,
\tP\in \cD_X[s,\rd_t^{-1}].
\leqno(4.1.2)
$$

Put $n=\dim X$. Then

\medskip\noindent
{\bf 4.2.~Theorem.}
$\tR_{f} \subset [\tal_{f},n-\tal_{f}],\quad
\tm_{\alpha} \le n-\tal_{f} - \alpha + 1$.

\medskip
(The proof uses the filtered duality for $\varphi_f$, see [35].)

\medskip\noindent
{\bf 4.3.~Spectrum.}  We define the spectrum by
$\Sp(f,x)=\sum_{\alpha}n_{\alpha}t^{\alpha}$ with
$$
n_{\alpha}:=\msum_j(-1)^{j-n+1}\dim
\Gr_F^p\tH^j(F_x,\bC)_{\lambda},
\leqno(4.3.1)
$$
where $p=[n-\alpha],\,\lambda=e^{-2\pi i\alpha}$, and $F$ is
the Hodge filtration (see [12]) of the mixed Hodge structure
on the Milnor cohomology,
see [44].
We define
$$
E_f=\{\alpha\mid n_{\alpha}\ne 0\}\,\,\, \text{(called the exponents)}.
\leqno(4.3.2)
$$

\medskip\noindent
{\bf 4.4.~Remarks.} (i)
If $f$ has an isolated singularity at the origin, then $\tal_{f,x}$
coincides with the minimal exponent as a corollary of results of
Malgrange [26], Varchenko [45], Scherk-Steenbrink [41].

\medskip
(ii) If $f$ is weighted-homogeneous with an isolated
singularity at the origin, then by Kashiwara (unpublished)
$$
\tR_f=E_f,\quad
\max \tR_f=n-\tal_f,\quad
\tm_\alpha=1\,\,(\alpha\in\tR_f).
\leqno(4.4.1)
$$
If $f=\sum_i x_i^2$, then $\tal_f=n/2$ and this follows
from the above Theorem~(4.2).

By Steenbrink [42], we have moreover
$$
\Sp(f,x)=\mprod_i (t-t^{w_i})/(t^{w_i}-1),
\leqno(4.4.2)
$$
where $(w_1,\dots,w_n)$ is the weights of $f$, i.e.
$f$ is a linear combination of monomials $x_1^{m_1}\cdots x_n^{m_n}$
with $\sum_i w_i m_i=1$.

\medskip\noindent
{\bf 4.5.~Malgrange's formula} (isolated singularities case).
We have the Brieskorn lattice [5] and its saturation defined by
$$
H''_f = \Omega_{X,x}^n/df\wedge d\Omega_{X,x}^{n-2},\quad
\tH''_f=\msum_{i\ge 0}(t\rd_t)^iH''_f \subset H''_f[t^{-1}].
\leqno(4.5.1)
$$
These are finite $\bC\{t\}$-modules with a regular singular connection.

\medskip\noindent
{\bf 4.6.~Theorem} (Malgrange [26]). {\it
The reduced $b$-function $\tb_f(s)$ coincides with the minimal
polynomial of $-\rd_tt$ on $\tH''_f/t\tH''_f$.
}

\medskip
(The above formula of Kashiwara on $b$-function (4.4.1) can be
proved by using this together with Brieskorn's calculation.)

\medskip\noindent
{\bf 4.7.~Asymptotic Hodge structure}
(Varchenko [45], Scherk-Steenbrink [41]). {\it
In the isolated singularity case we have
$$
F^pH^{n-1}(F_x,\bC)_{\lambda}=\Gr_V^{\alpha}H''_f,
\leqno(4.7.1)
$$
using the canonical isomorphism
$$
H^{n-1}(F_x,\bC)_{\lambda}=\Gr_V^{\alpha}H''_f[t^{-1}],
\leqno(4.7.2)
$$
where $p=[n-\alpha],\lambda=e^{-2\pi i\alpha}$,
and $V$ on $H''_f[t^{-1}]$ is
the filtration of Kashiwara and Malgrange.
}

\medskip
(This can be generalized to the non-isolated singularity case
using mixed Hodge modules.)

\medskip\noindent
{\bf 4.8.~Reformulation of Malgrange's formula.}
We define
$$
\tF^pH^{n-1}(F_x,\bC)_{\lambda}=\Gr_V^{\alpha}\tH''_f,
\leqno(4.8.1)
$$
using the canonical isomorphism (4.7.2),
where $p=[n-\alpha],\lambda=e^{-2\pi i\alpha}$.
Then
$$
\tm_{\alpha}=\hbox{
the minimal polynomial of $N$ on
$\Gr_{\tF}^pH^{n-1}(F_x,\bC)_{\lambda}$}.
\leqno(4.8.2)
$$

\medskip\noindent
{\bf 4.9.~Remark.}
If $f$ is weighted homogeneous with an isolated singularity, then
$$
\tF=F,\quad \tR_f = E_f\,\,\,\hbox{(by Kashiwara)}.
\leqno(4.9.1)
$$
If $f$ is not weighted homogeneous (but with isolated singularities),
then
$$
\tR_f\subset \mcup_{k\in\bN}(E_f-k),\,\,\,
\tal_f=\min \tR_f=\min E_f.
\leqno(4.9.2)
$$

\medskip\noindent
{\bf 4.10.~Example.} If $f=x^5+y^4+x^3y^2$, then
$$
E_f=\Bigl\{\frac{i}{5}+\frac{j}{4}: 1\le i \le 4,\,1\le j\le 3\Bigr\},
\quad
\tR_f=E_f\cup  \Big\{\frac{11}{20}\Big\}\setminus \Big\{\frac{31}{20}\Big\}.
$$

More generally, if $f=g+h$ with $g$ weighted homogeneous
and $h$ is a linear combination of monomials of higher degrees, then
$E_f=E_g$ but $\tR_f\ne\tR_g$ if $f$ is a non trivial
deformation.

\medskip\noindent
{\bf 4.11.~Relation with rational singularities} [34].
{\it Assume $D:=f^{-1}(0)$ is reduced.
Then $D$ has rational singularities if and only if $\tal_f>1$.
Moreover, $\omega_D/\rho_*\omega_{\tD}\simeq
F_{1-n}\varphi_f\cO_X$,
where $\rho:\tD\to D$ is a resolution of singularities.
}

\medskip
In the isolated singularities case, this was proved in 1981 (see [31])
using the coincidence of $\tal_f$ and the minimal exponent.

\medskip\noindent
{\bf 4.12.~Relation with the pole order filtration} [34].
{\it Let $P$ be the pole order filtration on $\cO_X(*D)$, i.e.
$P_i=\cO_X((i+1)D)$ if $i\ge 0$, and $P_i=0$ if $i<0$.
Let $F$ be the Hodge filtration on $\cO_X(*D)$.
Then $F_i\subset P_i$ in general, and
$F_i=P_i\,$ on a neighborhood of $x$ for $\,i\le\tal_{f,x} - 1$.
}

\medskip
(For the proof we need the theory of microlocal $b$-functions [35].)

\medskip\noindent
{\bf 4.13.~Remark.}  In case $X=\bP^n$, replacing $\tal_{f,x}$ with
$[(n-r)/d]$ where $r = \dim\Sing \,D$ and $d = \deg D$,
the assertion was obtained by Deligne (unpublished).

\bigskip\bigskip
\centerline{\bf 5. Relation with multiplier ideals.}

\bigskip\noindent
{\bf 5.1.~Multiplier ideals.}
Let  $D=f^{-1}(0)$, and $\cJ(X,\alpha D)$ be the multiplier
ideals for $\alpha\in\bQ$, i.e.
$$
\cJ(X,\alpha D)=\rho_*\omega_{\tX/X}(-\msum_i[\alpha m_i]\tD_i)),
\leqno(5.1.1)
$$
where $\rho:(\tX,\tD)\to(X,D)$ is an embedded resolution
and $\tD=\sum_i m_i\tD_i:=\rho^*D$.
There exist jumping numbers $0<\alpha_0<\alpha_1<\cdots$ such that
$$
\cJ(X,\alpha_j D)=\cJ(X,\alpha D)\ne\cJ(X,\alpha_{j+1} D)\quad
\text{for}\quad \alpha_j\le\alpha<\alpha_{j+1}.
\leqno(5.1.2)
$$

Let $V$ denote also the induced filtration on
$$
\cO_X\subset \cO_X[\rd_t]\delta(f-t).
$$

\medskip\noindent
{\bf 5.2.~Theorem} (Budur, S.\ [10]). {\it
If $\alpha$ is not a jumping number,
$$
\cJ(X,\alpha D)=V^{\alpha}\cO_X.
\leqno(5.2.1)
$$
For $\alpha$ general we have for $0<\varepsilon\ll 1$
}
$$
\cJ(X,\alpha D)=V^{\alpha+\varepsilon}\cO_X,\quad
V^{\alpha}\cO_X=\cJ(X,(\alpha-\varepsilon) D).
\leqno(5.2.2)
$$

\medskip
Note that $V$ is left-continuous and
$\cJ(X,\alpha D)$ is right-continuous, i.e.
$$
V^{\alpha}\cO_X=V^{\alpha-\varepsilon}\cO_X,\quad
\cJ(X,\alpha D)=\cJ(X,(\alpha+\varepsilon) D).
\leqno(5.2.3)
$$

The proof of (5.2) uses the theory of bifiltered direct images [32], [33]
to reduce the assertion to the normal crossing case.

As a corollary we get another proof of the results of Ein, Lazarsfeld,
Smith and Varolin [16], and of Lichtin, Yano and Koll\'ar [25]:

\bigskip
\noindent
{\bf 5.3.~Corollary.}

\smallskip\noindent
(i) $\{\hbox{Jumping numbers}\le 1\}\subset R_f$, see [16].

\smallskip\noindent
(ii) $\alpha_f=$ minimal jumping number, see [25].

\medskip
Define $\alpha'_{f,x}=\min_{y\ne x}\{\alpha_{f,y}\}$. Then

\medskip\noindent
{\bf 5.4.~Theorem.} {\it
If $\xi f =f$ for a vector field $\xi$, then
}
$$
R_f\cap(0,\alpha'_{f,x})=\{\hbox{Jumping numbers}\}\cap
(0,\alpha'_{f,x}).
\leqno(5.4.1)
$$

\medskip
(This does not hold without the assumption on
$\xi$ nor for $[\alpha'_{f,x},1)$.)

\medskip
For the constantness of the jumping numbers under a topologically
trivial deformation of divisors, see [14].

\bigskip\bigskip
\centerline{\bf 6. $b$-Functions for any subvarieties.}

\bigskip\noindent
{\bf 6.1.}
Let $Z$ be a closed subvariety of a smooth $X$, and
$f=(f_1,\dots,f_r)$ be generators of the ideal of $Z$
(which is not necessarily reduced nor irreducible).
Define the action of $ t_{j} $ on
$$
\cO_{X}\bigl[\hbox{$\frac{1}{f_1\cdots f_r}$}\bigr][s_{1},\dots, s_{r}]
\mprod_{i}f_{i}^{{s}_{i}},
$$
by
$ t_{j}(s_{i}) = s_{i} + 1 $ if
$ i = j $, and
$ t_{j}(s_{i}) = s_{i} $ otherwise.
Put $ s_{i,j} := s_{i}t_{i}^{-1}t_{j} $,
$s=\sum_i s_i$.
Then $b_f(s)$ is the monic polynomial of the least degree satisfying
$$
b_{f}(s)\mprod_{i}f_{i}^{{s}_{i}} =
\msum_{k=1}^r P_{k}t_{k}\mprod_{i}f_{i}^{{s}_{i}},
\leqno(6.1.1)
$$
where
$ P_{k} $ belong to the ring generated by
$ \cD_{X} $ and
$ s_{i,j} $.

\medskip
Here we can replace $\prod_{i}f_{i}^{{s}_{i}}$ with
$\prod_{i}\delta(t_i-f_i)$,
using the direct image by the graph of $f:X\to\bC^r$.
Then the existence of $b_f(s)$ follows from the theory
of the $V$-filtration of Kashiwara and Malgrange.
This $b$-function has appeared in work of Sabbah [30] and
Gyoja [18] for the study of $b$-functions of several variables.

\medskip\noindent
{\bf 6.2.~Theorem} (Budur, \Mustata, S.\ [8]). {\it
Let $c=\hbox{codim}_XZ$.
Then $b_Z(s):=b_f(s-c)$ depends only on $Z$ and is independent
of the choice of $f=(f_1,\dots,f_r)$ and also of $r$.
}

\medskip\noindent
{\bf 6.3.~Equivalent definition.}
The $b$-function
$ b_{f}(s) $ coincides with the monic polynomial of the least degree
satisfying
$$
b_{f}(s)\mprod_{i}f^{s_{i}} \in\msum_{|c|=1}
\cD_{X}[\bs]\,
\mprod_{c_i<0}\hbox{$\binom{s_{i}}{-c_{i}}$}
\mprod_{i}f_{i}^{s_{i}+c_{i}},
\leqno(6.3.1)
$$
where
$ c = (c_{1}, \dots, c_{r}) \in \bZ^{r} $ with
$ |c|:=\sum_{i}c_{i} = 1 $.
Here $\cD_{X}[\bs]=\cD_{X}[s_1,\cdots,s_r]$.

This is due to \Mustata, and is used in the monomial ideal case.
Note that the well-definedness does not hold
without the term $\prod_{c_i<0}\binom{s_{i}}{-c_{i}}$.

We have the induced filtration $V$ by
$$
\cO_X\subset i_{f+}\cO_X=\cO_X[\rd_1,\dots,\rd_r]\mprod_i
\delta(t_i-f_i).
$$

\medskip\noindent
{\bf 6.4.~Theorem} (Budur, \Mustata, S.\ [8]). {\it
If $\alpha$ is not a jumping number,
$$
\cJ(X,\alpha Z)=V^{\alpha}\cO_X.
\leqno(6.4.1)
$$
For $\alpha$ general we have for $0<\varepsilon\ll 1$
}
$$
\cJ(X,\alpha Z)=V^{\alpha+\varepsilon}\cO_X,\quad
V^{\alpha}\cO_X=\cJ(X,(\alpha-\varepsilon) Z).
\leqno(6.4.2)
$$

\medskip\noindent
{\bf 6.5.~Corollary} (Budur, \Mustata, S.\ [8]). {\it
We have the inclusion}
$$
\{\hbox{Jumping numbers}\}\cap[\alpha_f,\alpha_f+1)\subset R_f.
\leqno(6.5.1)
$$

\medskip\noindent
{\bf 6.6.~Theorem} (Budur, \Mustata, S.\ [8]). {\it
If $Z$ is reduced and is a local complete intersection, then
$Z$ has only rational singularities if and only if $\alpha_f=r$ with multiplicity $1$.
}

\bigskip\bigskip
\centerline{\bf 7. Monomial ideal case.}

\bigskip\noindent
{\bf 7.1.~Definition.} Let
$ \fa \subset\bC[x] := \bC[x_{1},\dots, x_{n}] $ a monomial ideal.
We have the associated semigroup defined by
$$
\Gamma_{\fa} = \{u\in\bN^n\mid x^u\in\fa\}.
$$
Let
$ P_{\fa} $ be the convex hull of $\Gamma_{\fa}$ in
$ \bR_{\ge 0}^{n} $.
For a face $Q$ of $P_{\fa}$, define

\medskip\quad
$M_{Q}$ : the subsemigroup of $\bZ^n$ generated by
$u-v$ with $u\in\Gamma_{\fa}$, $v\in\Gamma_{\fa}\cap Q$.

\medskip\quad
$M_Q'=v_0+M_Q$ for $v_0\in\Gamma_{\fa}\cap Q$
(this is independent of $v_0$).

\medskip
For a face $Q$ of $P_{\fa}$ not contained in any coordinate
hyperplane, take a linear function with rational coefficients
$L_Q:\bR^n\to\bR$ whose restriction to $Q$ is 1.
Let

\medskip\quad
$V_Q$ : the linear subspace generated by $Q$.

\smallskip\quad
$e=(1,\dots,1)$.

\smallskip\quad
$R_Q=\{L_Q(u)\mid u\in (e+(M_Q\setminus M_{Q}'))\cap V_Q\}$,

\smallskip\quad
$R_{\fa}=\{$roots of $b_{\fa}(-s)\}$.

\medskip\noindent
{\bf 7.2.~Theorem} (Budur, \Mustata, S.\ [9]). {\it We have
$R_{\fa}=\bigcup_QR_Q$ with $Q$ faces of $P_{\fa}$
not contained in any coordinate hyperplanes.
}

\medskip\noindent
{\it Outline of the proof.}
Let $f_j=\prod_ix_i^{a_{i,j}}$, $\ell_i(\bs)=\sum_j a_{i,j}s_j$.
Define
$$
\hbox{$g_c(\bs)=\mprod_{c_i<0}\binom{s_{i}}{-c_{i}}
\mprod_{\ell_i(c)>0}\binom{\ell_i(\bs)+\ell_i(c)}{\ell_i(c)}$}.
$$
Let
$I_{\fa}\subset \bC[\bs]$ be the ideal generated by $g_c(\bs)$ with
$c\in\bZ^r,\sum_i c_i=1$.
Then

\medskip\noindent
{\bf 7.3.~Proposition} (\Mustata). {\it
The $b$-function
$b_{\fa}[s]$ of the monomial ideal $\fa$
is the monic generator of $\bC[s]\cap I_{\fa}$,
where $s=\sum_i s_i$.
}

\medskip
Using this, Theorem (7.2) follows from elementary computations.

\medskip\noindent
{\bf 7.4.~Case $n=2$.} 
Here it is enough to consider only 1-dimensional $Q$ by (7.2).
Let
$Q$ be a compact face of $P_{\fa}$ with $\{v^{(1)},v^{(2}\}=\rd Q$,
where $v^{(i)}=(v^{(i)}_1,v^{(i)}_2)$ with $v^{(1)}_1<v^{(2)}_1$,
$v^{(1)}_2>v^{(2)}_2$.
Let

\medskip\quad
$G$ : the subgroup generated by $u-v$ with $u,v\in Q\cap \Gamma_{\fa}$.

\smallskip\quad
$v^{(3)}\in Q\cap\bN^2$ such that $v^{(3)}-v^{(1)}$ generates $G$.

\smallskip\quad
$S_Q=\{(i,j)\in\bN^2\mid i< v^{(3)}_1, \,j< v^{(1)}_2\}$.

\smallskip\quad
$S^{[1]}_Q=S\cap M_{Q}'$, \,\,$S^{[0]}_Q=S_Q\setminus S^{[1]}_Q$.

\medskip\noindent
Then

\medskip\quad
$R_Q=\{L_Q(u+e)-k\mid u\in S^{[k]}_Q\,(k=0,1)\}.$

\medskip\noindent
In the case $Q\subset\{x=m\}$, we have $R_Q=\{i/m\mid i=1,\dots,m\}$.

\medskip\noindent
{\bf 7.5.~Examples.} (i) If $\fa=(x^ay,xy^b)$, with
$ a,b \ge 2 $, then
$$
R_{\fa}=\Bigl\{{\frac{(b - 1)i + (a - 1)j}{ab - 1}}\,\Big|\,
1\le i\le a,\,\, 1\le j\le b\Bigr\}.
$$

\medskip\noindent
(ii) If $\fa=(xy^5,x^3y^2,x^5y)$, then
$S_Q^{[1]}=\emptyset$ and
$$
R_{\fa}=\Bigl\{\frac{5}{13},\,\frac{i}{13}\,(7\le i\le 17),\,
\frac{19}{13},\,\frac{j}{6}\,\,(3\le j\le 9)\Bigr\}.
$$

\medskip\noindent
(iii) If $\fa=(xy^5,x^3y^2,x^4y)$, then
$S_Q^{[1]}=\{(2,4)\}$ for $\rd Q=\{(1,5),(3,2)\}$ with
$L_Q(v_1,v_2)=(3v_1+2v_2)/13$, and
$$
R_{\fa}=\Bigl\{\frac{i}{13}\,\,(5\le i\le 17),\,\,\frac{j}{5}\,\,
(2\le j\le 6)\Bigr\}.
$$
Here $19/13$ is shifted to $6/13$.

\medskip\noindent
{\bf 7.6.~Comparison with exponents.}
If $n=2$ and $f$ has a nondegenerate Newton polygon
with compact faces $Q$, then by Steenbrink [43]
$$
E_f\cap(0,1]=\mcup_Q E_Q^{\le 1}\quad\text{with}\quad
E_Q^{\le 1}=\{L_Q(u)\mid u\in\overline{\{0\}\cup Q}\cap \bZ_{>0}^2\},
$$
where $\overline{\{0\}\cup Q}$ is the convex hull of
$\{0\}\cup Q$.
Here we have the symmetry of $E_f$ with center 1.

\medskip\noindent
{\bf 7.7.~Another comparison.}
If $\fa=(x_1^{a_1},\dots,x_n^{a_n})$, then
$$
R_{\fa}=\bigl\{\msum_i p_i/a_i\mid 1\le p_i\le a_i\bigr\}.
$$
On the other hand, if $f=\sum_i x_i^{a_i}$, then
$$
\tR_f=E_f=\bigl\{\msum_i p_i/a_i\mid 1\le p_i\le a_i-1\bigr\}.
$$

\bigskip\bigskip
\centerline{\bf 8. Hyperplane arrangements.}

\bigskip\noindent
{\bf 8.1.} Let
$D$ be a central hyperplane arrangement in $X=\bC^n$.
Here, central means an affine cone of $Z\subset\bP^{n-1}$.
Let $f$ be the reduced equation of $D$ and $d:=\deg f>n$.
Assume $D$ is not the pull-back of $D'\subset\bC^{n'}(n'<n)$.

\medskip\noindent
{\bf 8.2.~Theorem.} (i) $\max R_f<2-\frac{1}{d}$. (ii) $m_1=n$.

\medskip
Proof of (i) uses a partial generalization of a solution of
Aomoto's conjecture due to Esnault, Schechtman, Viehweg, Terao,
Varchenko ([17], [40]) together with a generalization of Malgrange's
formula (4.8) as below:

\medskip\noindent
{\bf 8.3.~Theorem} ({\rm Generalization of Malgrange's formula}) [36].
{\it There exists a pole order filtration $P$ on
$H^{n-1}(F_0,\bC)_{\lambda}$ such that if $(\alpha+\bN)\cap
R'_f=\emptyset$, then
$$
\alpha\in R_f\Leftrightarrow \Gr_P^pH^{n-1}(F_0,\bC)_{\lambda}\ne 0,
\leqno(8.3.1)
$$
with $p=[n-\alpha],\lambda=e^{-2\pi i\alpha}$,
where $R'_f=\cup_{x\ne 0}R_{f,x}$.
}

\medskip
This reduces the proof of (8.2)(i) to
$$
P^iH^{n-1}(F_0,\bC)_{\lambda}=H^{n-1}(F_0,\bC)_{\lambda},
\leqno(8.3.2)
$$
for $i=n-1$ if $\lambda=1$ or $e^{2\pi i/d}$,
and $i=n-2$ otherwise.

\medskip\noindent
{\bf 8.4.~Construction of the pole order filtration P.}
Let $U=\bP^{n-1}\setminus Z$, and $F_0=f^{-1}(0)\subset\bC^n$.
Then $F_0=\tU$ with $\pi:\tU\to U$ a $d$-fold covering ramified over $Z$.
Let $L^{(k)}$ be the local systems of rank 1 on $U$ such that
$\pi_*\bC=\moplus_{0\le i<d}L^{(k)}$
and $T$ acts on $L^{(k)}$ by $e^{-2\pi i k/d}$. Then
$$
H^j(U,L^{(k)})=H^j(F_0,\bC)_{{\bf e}(-k/d)},
\leqno(8.4.1)
$$
and $P$ is induced by the pole order filtration on
the meromorphic extension $\cL^{(k)}$ of $L^{(k)}\otimes_{\bC}\cO_U$
over $\bP^{n-1}$, see [15], [36], [37].
This is closely related to:

\medskip\noindent
{\bf 8.5.~Solution of Aomoto's conjecture} ([17], [40]). Let
$Z_i$ be the irreducible components of $Z\,\,(1\le i\le d)$,
$g_i$ be the defining equation of $Z_i$ on $\bP^{n-1}\setminus Z_d\,\,
(i<d)$, and
$\omega:=\sum_{i<d}\alpha_i\omega_i$ with $\omega_i=dg_i/g_i$,
$\alpha_i\in\bC$.
Let $\nabla$ be the connection on $\cO_U$ such that
$\nabla u=du+\omega\wdg u$.
Set $\alpha_d=-\sum_{i<d}\alpha_i$.
Then $H_{\rm DR}^{\ssbull}(U,(\cO_U,\nabla))$ is
calculated by
$$
(\cA^{\ssbull}_{\alpha},\omega\wedge) \quad\text{with}\quad
\cA^p_{\alpha}=\msum \bC\omega_{i_1}\wdg\cdots\wdg\omega_{i_p},
$$
if $\sum_{Z_i\supset L}\alpha_i\notin \bN\setminus\{0\}$ for any
{\it dense} edge $L\subset Z$ (see (8.7) below).
Here an edge is an intersection of $Z_i$.

\medskip
For the proof of (8.2)(ii) we have

\medskip\noindent
{\bf 8.6.~Proposition.} {\it $N^{n-1}\psi_{f,\lambda}\bC\ne 0$ if
$\Gr^W_{2n-2}H^{n-1}(F_x,\bC)_{\lambda}\ne 0$.}

\medskip
(Indeed,
$N^{n-1}:\Gr^W_{2n-2}\psi_{f,\lambda}\bC\simto
\Gr^W_0\psi_{f,\lambda}\bC$ by the definition of $W$, and
the assumption of (8.6) implies
$\Gr^W_{2n-2}\psi_{f,\lambda}\bC\ne 0$.)

\medskip
Then we get (8.2)(ii), since
$\omega_{i_1}\wdg\cdots\wdg\omega_{i_{n-1}}\ne 0$ in
$\Gr^W_{2n-2}H^{n-1}(\bP^{n-1}\setminus Z,\bC)=
\Gr^W_{2n-2}H^{n-1}(F_x,\bC)_{1}$.

\medskip\noindent
{\bf 8.7.~Dense edges.}
Let $D=\cup_iD_i$ be the irreducible decomposition.
Then $L=\cap_{i\in I}D_i$ is called an edge of
$D\,\,\,(I\ne \emptyset$),

We say that an edge $L$ is {\it dense} if
$\{D_i/L\,|\,D_i\supset L\}$ is indecomposable.
Here $\bC^n\supset D$ is called decomposable if
$\bC^n=\bC^{n'}\times\bC^{n''}$ such that
$D$ is the union of the pull-backs from $\bC^{n'},
\bC^{n''}$ with $n',n''\ne 0$.

Set $m_L=\#\{D_i\,|\,D_i\supset L\}$. For $\lambda\in\bC$, define
$$
\cD\cE(D) = \{\hbox{dense edges of }D\},\quad
\cD\cE(D,\lambda)=\{L\in\cD\cE(D)\,|\,\lambda^{m_L}=1\}.
$$
We say that $L,\,L'$ are {\it strongly adjacent} if
$L\subset L'$ or $L\supset L'$ or $L\cap L'$ is non-dense.
Let
$$
\aligned
m(\lambda)&=\max\{|S|\,\big|\,
S\subset\cD\cE(D,\lambda) \,\, \text{such that}
\\
&\qquad\text{any $L,L'\in S$ are strongly adjacent}\}.
\endaligned
$$

\medskip\noindent
{\bf 8.8.~Theorem} [37]. {\it $m_{\alpha}\le m(\lambda)$ with
$\lambda=e^{-2\pi i\alpha}$.}

\medskip\noindent
{\bf 8.9.~Corollary.}
$R_f\subset\bigcup_{L\in\cD\cE(D)}\bZ m_L^{-1}$.

\medskip\noindent
{\bf 8.10.~Corollary.} {\it If {\rm GCD}$(m_L,m_{L'})=1$ for any
strongly adjacent $L,L'\in\cD\cE(D)$, then $m_{\alpha}=1$ for any
$\alpha\in R_f\setminus \bZ$.}

\medskip
Theorem 2 follows from the canonical resolution of singularities
$\pi:(\tX,\tD)\to (\bP^{n-1},D)$ due to [40],
which is obtained by blowing up along
the proper transforms of the dense edges.
Indeed, mult $\tD(\lambda)_{\rm red}\le m(\lambda)$,
where $\tD(\lambda)$ is the union of $\tD_i$ such that
$\lambda^{\tm_i}=1$ and $\tm_i=\hbox{mult}_{\tD_i}\tD$.

\medskip\noindent
{\bf 8.11.~Theorem} (\Mustata\,[29]). For a central arrangement,
$$
\hbox{$\cJ(X,\alpha D)=I_0^k$ with $k=[d\alpha]-n+1$ if
$\alpha<\alpha'_f$},
\leqno(8.11.1)
$$
where $I_0$ is the ideal of $0$ and
$\alpha'_f=\min_{x\ne 0}\{\alpha_{f,x}\}$.

\medskip
(This holds for the affine cone of any divisor on $\bP^{n-1}$,
see [36].)

\medskip\noindent
{\bf 8.12.~Corollary.} {\it We have
$\dim F^{n-1}H^{n-1}(F_0,\bC)_{{\bf e}(-k/d)}=
\binom{{k-1}}{{n-1}}$
for $0<\frac{k}{d}<\alpha'_f$,
and the same holds with $F$ replaced by $P$.
}

\medskip\noindent
{\bf 8.13.~Corollary.}
$\alpha_{f} =\min(\alpha'_f,\frac{n}{d})<1$.

\medskip
(Note that $\alpha_f$ coincides with the minimal jumping number.)

\medskip\noindent
{\bf 8.14.~Generic case.} If $D$ is a generic central
hyperplane arrangement, then
$$
b_f(s)=(s+1)^{n-1}\mprod_{j=n}^{2d-2}(s+\hbox{$\frac{j}{d}$})
\leqno(8.14.1)
$$
by U.~Walther [46] (except for the multiplicity of $-1$).
He uses a completely different method.

Note that Theorems (8.2) and (8.8) imply that the left-hand side
divides the right-hand side of (8.14.1), and the equality follows
using also (8.12).

\medskip\noindent
{\bf 8.15.~Explicit calculation.}
Let $\alpha=k/d$,
$\lambda=e^{-2\pi i\alpha}$ for $k\in \{1,\dots,d\}$.
If $\alpha\ge \alpha'_f$, we assume there is $I\subset\{1,\dots,d-1\}$
such that $|I|=k-1$, and the condition of [40]
$$
\hbox{$\sum_{Z_i\supset L}\alpha_i\notin \bN\setminus\{0\}$ for any
dense edge $L\subset Z$,}
\leqno(8.15.1)
$$
is satisfied for
$$
\hbox{$\alpha_i=1-\alpha$ if $i\in I\cup\{d\}$, and $-\alpha$
otherwise}.
\leqno(8.15.2)
$$

Let $V(I)$ be the subspace of $H^{n-1}\cA^{\ssbull}_{\alpha}$
generated by
$$
\omega_{i_1}\wdg\cdots\wdg\omega_{i_{n-1}}\quad\text{for}\quad
\{i_1,\dots,i_{n-1}\}\subset I.
$$

\medskip\noindent
{\bf 8.16.~Theorem.} {\it Let $\alpha=k/d$,
$\lambda=e^{-2\pi i\alpha}$ for $k\in \{1,\dots,d\}$. Then

\noindent
{\rm (a)} If $k=d-1$ or $d$, then
$\alpha\in R_f$, $\alpha+1\notin R_f$.

\noindent
{\rm (b)} If $\alpha<\alpha'_f$, then $\alpha\in R_f\Leftrightarrow
k\ge d$.

\noindent
{\rm (c)} If $\binom{k-1}{n-1}<\dim H^{n-1}(F_0,\bC)_{\lambda}$, then
$\alpha+1\in R_f$.

\noindent
{\rm (d)} If $\alpha<\alpha'_f$, $\alpha\notin R'_f+\bZ$ and
$\binom{k-1}{n-1}=\chi(U)$, then $\alpha+1\notin R_f$.

\noindent
{\rm (e)} If $\alpha\ge \alpha'_f$ and $V(I)\ne 0$, then $\alpha\in R_f$.

\noindent
{\rm (f)} If $\alpha\ge \alpha'_f$ and $V(I)=H^{n-1}
\cA^{\ssbull}_{\alpha}$, then $\alpha+1\notin R_f$.
}

\medskip\noindent
{\bf 8.17.~Theorem} [37]. {\it 
Assume $n=3$, mult$_zZ\le 3$ for any $z\in Z\subset \bP^2$, and
$d\le 7$. Let $\nu_3$ be the number of triple points of $Z$, and
assume $\nu_3\ne 0$. Then
$$
\hbox{$b_f(s)=(s+1)\prod_{i=2}^4(s+\frac{i}{3})\,
\prod_{j=3}^{r}(s+\frac{j}{d})$},
\leqno(8.17.1)
$$
with $r=2d-2$ or $2d-3$.
We have $r=2d-2$ if $\nu_3<d-3$, and the converse holds
for $d<7$. In case $d=7$, we have $r=2d-3$ for $\nu_3>4$, however,
for $\nu_3=4$, $r$ can be both $2d-2$ and $2d-3$.
}

\medskip\noindent
{\bf 8.18.~Remarks.} (i) We have $\nu_3 < d-3$ if and only if
$$
\chi(U)=\hbox{$\frac{(d-2)(d-3)}{2}$}-\nu_3>
\hbox{$\frac{(d-3)(d-4)}{2}=\binom{d-3}{2}$}.
\leqno(8.18.1)
$$

\medskip\noindent
(ii) By (8.4.1) we have
$\chi(U)=h^2(F_0,\bC)_{\lambda}-h^1(F_0,\bC)_{\lambda}$
if $\lambda^d=1$ and $\lambda\ne 0$.

\medskip\noindent
(iii) Let $\nu'_i$ be the number of $i$-ple points of
$Z':=Z\cap\bC^2$. Then by [6]
$$
\hbox{$b_0(U)=1,\quad b_1(U)=d-1,\quad b_2(U)=\nu'_2+2\nu'_3,$}
\leqno(8.18.2)
$$

\medskip\noindent
{\bf 8.19.~Examples.} (i) For $(x^2-1)(y^2-1)=0$ in $\bC^2$ with $d=5$,
(8.17.1) holds with $r=7$, and $8/5\notin R_f$.
In this case we do not need to take $I$, because
$(d-2)/d=3/5<\alpha'_f=2/3$.
We have $b_1(U)=b_2(U)=4$ and
$h^2(F_0,\bC)_{\lambda}=\chi(U)=1$ if
$\lambda^5=1$ and $\lambda\ne 1$. 
So $j/5\in R_f$ for $3\le j\le 7$ by (a), (b), (c),
and $8/5\notin R_f$ by (d).

\medskip\noindent
(ii) For $(x^2-1)(y^2-1)(x+y)=0$ in $\bC^2$ with $d=6$,
(8.17.1) holds with $r=9$, and $10/6\notin R_f$.
In this case we have
$b_1(U)=5,b_2(U)=6,\chi(U)=2,
h^1(F_0,\bC)_{\lambda}=1, h^2(F_0,\bC)_{\lambda}=3$ for
$\lambda=e^{\pm 2\pi i/3}$.
Then $4/6\in R_f$ by (e) and $10/6\notin R_f$ by (f),
where $I^c$ corresponds to $(x+1)(y+1)=0$.
For other $j/6$, the argument is the same as in (i).

\medskip\noindent
(iii) For $(x^2-y^2)(x^2-1)(y+2)=0$ in $\bC^2$ with $d=6$,
(8.17.1) holds with $r=10$, and $10/6\in R_f$.
In this case we have
$b_1(U)=5, b_2(U)=9, \chi(U)=5,
h^1(F_0,\bC)_{\lambda}=0, \,\,\, h^2(F_0,\bC)_{\lambda}=5$ for
$\lambda=e^{\pm 2\pi i/3}$.
Then $4/6\in R_f$ by (e) and $10/6\in R_f$ by (c),
where $I^c$ corresponds to $(x+1)(y+2)=0$.

\medskip\noindent
(iv) For $(x^2-y^2)(x^2-1)(y^2-1)=0$ in $\bC^2$ with $d=7$,
(8.17.1) holds with $r=11$, and $12/7\notin R_f$.
In this case we have
$b_1(U)=6, b_2(U)=9, \chi(U)=4,h^2(F_0,\bC)_{\lambda}=4$ if
$\lambda^7=1\quad\text{and}\quad\lambda\ne 1$.
Then $5/7\in R_f$ by (e) and $12/7\notin R_f$ by (f),
where $I^c$ corresponds to $(x+1)(y+1)=0$.
Note that $5/7$ is not a jumping number.

\end{document}